\documentclass{amsart}
\usepackage{amssymb, amsfonts, lacromay}
\usepackage{mathrsfs,comment}
\usepackage[usenames,dvipsnames]{color}
\usepackage[normalem]{ulem}
\usepackage{url}

\title{A brief note on the Hopf and Kervaire invariant one problems}

\author{J.P. May}
\email{may@math.uchicago.edu}
\address{Department of Mathematics,
The University of Chicago, 
Chicago, IL 60637 USA}

\date{\today}

\begin{document}

\begin{abstract}  This  expository note advertises a little known but very old one-line solution of the Hopf invariant one problem and speculates on the (in)applicability of the idea to the Kervaire invariant one problem.
\end{abstract}

\maketitle

In 1958, Adams \cite{1958} introduced the Adams spectral sequence (ASS) in a failed attempt to solve the Hopf invariant one problem.   He observed that  $d_2(h_4) = h_0 h_3^2$ in the ASS since otherwise   $h_0h_3^2$ would survive to prove that $2 \sigma^2 \neq 0$ in the stable homotopy groups of spheres, a contradiction since $\sigma$ has odd degree.   Later in 1958 \cite{1958a}, with details published in 1960 \cite{1960}, he announced a solution of the Hopf invariant one problem that was obtained without use of the ASS.  He knew \cite[Lemma 2.3]{1958} early on that  to solve the problem it would suffice to prove that $d_2(h_n) = h_0 h_{n-1}^2$
for $n\geq 4$ in the ASS.  He overlooked a clever one-line inductive argument that John S. P. Wang \cite{Wang} found in 1967.  Together with a bit of calculation, it just  uses that
$$  d_2 (h_nh_{n+1}) = 0  \ \ \text{since} \ \    h_nh_{n+1} = 0.$$
The induction starts with $n=4$.  By the induction hypothesis,
$$  d_2(h_n) h_{n+1} = h_0h_{n-1}^2 h_{n+1}.$$
By that bit of calculation, very easy here in view of the low homological dimension $s=4$, that is non-zero.  Using \cite{MayThesis}, an easier proof than Wang's was available by 1967.  Therefore, by the Leibnitz rule, we must also have
$$h_n d_2(h_{n+1}) = h_0h_{n-1}^2 h_{n+1}. $$
This is non-zero and by a very easy inspection the only way that can be true is if
$$  d_2(h_{n+1})  = h_0 h_n^2.$$
This completes the induction.  For consistency, we note that $h_n^3 = h_{n-1}^2 h_{n+1}$.
That is all there is to the proof of Hopf invariant one!

Speculating, let us try to apply the same argument towards the Kervaire invariant one  problem, which concerns differentials on $h_n^2$ in the ASS.  We first quickly review what is known, bringing up to date the summary in \cite{Xu}.  If $h_n^2$ survives to $E_{\infty}$, we write $\theta_n$ for any element of $\pi_{2^{n+1}-2}S^0$ that detects it.  For $n=1$, $2$, and $3$, the unique elements $\theta_n$ are given by  $\eta^2$, $\nu^2$, and $\sigma^2$.  

Mahowald and Tangora \cite{MaTa} showed that $\theta_4$ exists.  In \cite{BJM2},  Barratt, Jones, and Mahowald  showed that $\theta_5$ exists, and in \cite{BJM1} they showed  that if there exists a $\theta_n$ such that $2\theta_n = 0$ and 
$\theta_n^2 = 0$, then there exists  a $\theta_{n+1}$ such that $2\theta_{n+1} =0$.  The element $\theta_4$ is unique, and the question of whether or not $\theta_4^2=0$ has a long and checkered history.  Xu \cite{Xu} proved that $\theta_4^2=0$, and he used that to give a new and simpler proof that $\theta_5$  exists and can be chosen so that $2\theta_5 = 0$.  Thus if  
$\theta_5^2=0$, then there exists a $\theta_6$ such that $2\theta_6 = 0$. However, it is not known that $\theta_5^2 = 0$.

Recently Weinan Lin, Guozheng Wang, and Zhouli Xu proved that $\tha_6$ exists  \cite{LWX}.  They use  highly nontrivial computer calculation to identify all possible non-zero diferentials  $d_r(h_6^2)$ and then proved by sheer ingenuity that none of them happen. 

By the earlier great work of Hill,  Hopkins, and Ravenel \cite{HHR}, for $n\geq 7$ all $h_n^2$ support some non-zero $d_r$, hence no $\theta_n$ exists.   It is not yet known for which $r$  $d_r(h_n^2) \neq 0$, let alone what that non-zero differential is. In view of the trivial solution of the Hopf invariant one problem, it is natural to dream that there might be a more informative inductive solution of the Kervaire invariant one problem.

We know that $d_r(h_7^2) \neq 0$ for some as yet unknown $r$ and it seems likely and provable that all  $h_n^2$  for $n>7$ survive to that same $E_r$.    Write  $x_n = d_r(h_n^2)$ for $n\geq 7$.  Assume inductively that $x_n \neq 0$.   The induction starts with $n=7$. Then 
$$  d_r (h_n^2h_{n+1}^2) = 0 \ \ \text{since} \ \   h_n^2h_{n+1}^2 = 0 .$$
By the induction hypothesis, 
$$d_r(h_n^2) h_{n+1}^2 = x_n h_{n+1}^2.$$
The hard  calculational part is to show that  $x_n h_{n+1}^2 \neq 0$.   Granting that, the Leibnitz rule gives
$$ h_n^2 d_r(h_{n+1}^2) = h_n^2 x_{n+1} = x_n h_{n+1}^2 \neq 0.$$
To complete the induction, the hard part then is to show that $ x_{n+1} h_{n+2}^2 \neq 0$.
Knowledge of $x_n$ should lead inductively to a description of $x_{n+1}$ in terms of Steenrod operations and/or Massey products involving $x_n$, and that could well lead to an inductive proof of the hard part.   This may seem to be out of reach calculationally but, not so long ago, so did determining whether or not $h_6^2$ survives.  

\bibliographystyle{plain}
\bibliography{references}

\end{document}